\documentclass[letterpaper,11pt]{article}%
\usepackage[textwidth=5.5in,textheight=9in,centering,letterpaper]{geometry}
\usepackage{amsxtra,amscd}
\usepackage[all,cmtip]{xy}
\usepackage{graphicx}
\usepackage[amsmath,hyperref,thmmarks]{ntheorem}
\usepackage{natbib}
\usepackage{verbatim}
\usepackage{paralist}
\usepackage{makeidx}
\usepackage{amsmath}
\usepackage{amsfonts}
\usepackage{amssymb}
\usepackage{times}
\theoremnumbering{arabic}
\theoremheaderfont{\scshape}
\RequirePackage{latexsym}
\theorembodyfont{\slshape}
\theoremseparator{}
\newtheorem{X}{X}[section]

\newtheorem{corollary}[X]{Corollary}

\newtheorem{lemma}[X]{Lemma}
\newtheorem{proposition}[X]{Proposition}
\newtheorem{theorem}[X]{Theorem}
\theorembodyfont{\upshape}

\newtheorem{definition}[X]{Definition}
\newtheorem{example}[X]{Example}

\newtheorem{question}[X]{Question}
\newtheorem{remark}[X]{Remark}
\newtheorem{plain}[X]{}

\theorembodyfont{\small}
\newtheorem*{nt}{Notes}
\theorembodyfont{\normalsize}

\theoremstyle{nonumberplain}
\theoremheaderfont{\sc}
\theorembodyfont{\normalfont}
\theoremsymbol{\ensuremath{_\Box}}
\RequirePackage{amssymb}
\newtheorem{proof}{Proof.}
\qedsymbol{\ensuremath{_\Box}}
\theoremclass{LaTeX}
\makeindex

\newcommand{\tsize}{\textstyle}

\newcommand\bquote{\begin{quote}}
\newcommand\equote{\end{quote}}
\newcommand\bsmall{\begin{small}}
\newcommand\esmall{\end{small}}

\let\cite=\citealt

\setdefaultitem{$\circ$\hspace{0.07in}}{}{}{}

\newcommand{\tstyle}{\textstyle}

\newcommand{\eb}[1]{{\itshape\bfseries#1}{\index{#1}}}
\renewcommand{\emph}{\eb}

\renewcommand{\bar}{\overline}
\renewcommand{\hat}{\widehat}

\renewcommand{\Gamma}{\varGamma}
\renewcommand{\Pi}{\varPi}
\renewcommand{\Sigma}{\varSigma}

\def\1{{1\mkern-7mu1}}

\DeclareMathOperator{\Br}{Br}

\DeclareMathOperator{\End}{End}

\DeclareMathOperator{\Gal}{Gal}

\DeclareMathOperator{\Hom}{Hom}

\DeclareMathOperator{\inv}{inv}

\DeclareMathOperator{\Nm}{Nm}

\DeclareMathOperator{\ord}{ord}

\DeclareMathOperator{\rank}{rank}

\newcommand{\Art}{\mathsf{Art}}

\newcommand{\Mot}{\mathsf{Mot}}

\begin{document}

\title{Motives over $\mathbb{F}_{p}$}
\author{J.S. Milne}
\date{July 22, 2006}
\maketitle

\begin{abstract}
In April, 2006, Kontsevich asked me whether the category of motives over
$\mathbb{F}_{p}$ ($p$ prime) has a fibre functor over a number field of
finite degree since he had a conjecture that more-or-less implied this. This
article is my response. Unfortunately, since the results are generally
negative or inconclusive, they are of little interest except perhaps for the
question they raise on the existence of a cyclic extension of $\mathbb{Q}$
having certain properties (see Question \ref{20e}).

\end{abstract}

Let $k$ be a finite field. Starting from any suitable class $\mathcal{S}{}$ of
algebraic varieties over $k$ including the abelian varieties and using the
correspondences defined by algebraic cycles modulo numerical equivalence, we
obtain a graded tannakian category $\Mot(k)$ of motives. Let $\Mot_{\text{0}%
}(k)$ be the subcategory of motives of weight $0$ and assume that the Tate
conjecture holds for the varieties in $\mathcal{S}{}$.

For a simple motive $X$, $D=\End(X)$ is a division algebra with centre the
subfield $F=\mathbb{Q}{}[\pi_{X}]$ generated by the Frobenius endomorphism
$\pi_{X}$ of $X$ and
\[
\rank(X)=[D\colon F]^{\frac{1}{2}}\cdot\lbrack F\colon\mathbb{Q}{}].
\]
Therefore, $D$ can act on a $\mathbb{Q}{}$-vector space of dimension
$\rank(X)$ only if it is commutative. Since this is never the case for the
motive of a supersingular elliptic curve or of the abelian variety obtained by
restriction of scalars from such a curve, there cannot be a $\mathbb{Q}{}%
$-valued fibre functor on the full category $\Mot(k)$. Let $k=\mathbb{F}{}%
_{q}$. Then, for each prime $v$ of $F$,%
\begin{equation}
\inv_{v}(D)=\left\{
\begin{array}
[c]{ll}%
1/2 & \text{if }v\text{ is real and }X\text{ has odd weight}\\
\dfrac{\ord_{v}(\pi_{X})}{\ord_{v}(q)}\cdot\lbrack F_{v}:\mathbb{Q}_{p}] &
\text{if }v|p\\
0 & \text{otherwise}%
\end{array}
\right.  \label{e6}%
\end{equation}
(Tate's formula; see \cite{milne1994}, 2.16). When $q=p$, $\ord_{v}(p)$ is the
ramification index $e(v/p)$, which divides the local degree $[F_{v}%
\colon\mathbb{Q}_{p}]$. Thus, for $k=\mathbb{F}{}_{p}$ and $X$ a motive of
weight $0$ (modulo $2)$, $D$ is commutative, and so the endomorphism algebras
provide no obstruction to $\Mot_{0}(\mathbb{F}{}_{p})$ being neutral. In this
note, we examine whether it is, in fact, neutral.

Before stating our results, we need some notations. Let $K$ be a CM subfield
of $\mathbb{C}$, finite and galois over $\mathbb{Q}$, and let $n$ be a
sufficiently divisible positive integer. Define $W^{K}(p,n)$ to be the group
of algebraic numbers $\pi$ in $\mathbb{C}{}$ such that

\begin{itemize}
\item $|\pi^{\prime}|=1$ for all conjugates $\pi^{\prime}$ of $\pi$ in
$\mathbb{C}{}$;

\item $p^{N}\pi$ is an algebraic integer for some $N$;

\item $\pi^{n}\in K$, and for every $p$-adic prime of $w$ of $K$,
$\dfrac{\ord_{w}(\pi^{n})}{n\cdot\ord_{w}(p)}[K_{w}\colon\mathbb{Q}{}_{p}%
]\in\mathbb{Z}$.
\end{itemize}

\noindent Define $\Mot_{0}^{K}(\mathbb{F}{}_{p},n)$ to be the category of
motives over $\mathbb{F}{}_{p}$ whose Weil numbers lie in $W^{K}(p,n)$. Let
$m$ be the order of $\mu(K)$. We prove the following.

\begin{description}
\item[(\ref{11a})] There exists a $\mathbb{Q}_{l}$-valued fibre functor on
$\Mot_{0}^{K}(\mathbb{F}{}_{p},n)$ for every prime $l$ of $\mathbb{Q}$
(including $p$ and $\infty$).

\item[(\ref{17})] There exists a $\mathbb{Q}$-valued fibre functor on
$\Mot_{0}^{K}(\mathbb{F}{}_{p},n)$ if and only if there exists a cyclic field
extension $L$ of $\mathbb{Q}$ of degree $mn$ such that

\begin{enumerate}
\item $(p)$ remains prime\footnote{By this I mean that the ideal generated by
$p$ in $\mathcal{O}{}_{L}$ is prime.} in $L$;

\item $\pi$ is a local norm at every prime $v$ of $\mathbb{Q}{}[\pi^{mn}]$
that ramifies in $\mathbb{Q}{}(\pi^{mn})\otimes_{\mathbb{Q}{}}L$.
\end{enumerate}

\noindent Moreover, we show that the generalized Riemann hypothesis sometimes
implies that there exists such an $L$.
\end{description}

Now consider the full category $\Mot_{0}(\mathbb{F}{}_{p})$ of motives of
weight $0$ over $\mathbb{F}{}_{p}$. Then
\[
\Mot_{0}(\mathbb{F}{}_{p})=\bigcup\nolimits_{K,n}\Mot_{0}^{K}(\mathbb{F}{}%
_{p},n),
\]
but the existence of a $\mathbb{Q}{}$-valued fibre functor on each of the
categories $\Mot_{0}^{K}(\mathbb{F}{}_{p},n)$ does not imply that there exists
a $\mathbb{Q}{}$-valued fibre functor on $\Mot_{0}(\mathbb{F}{}_{p})$. In
fact, we give a heuristic argument (due to Kontsevich) to show that there does
not exist such a fibre functor.

\textit{Throughout the article, we fix a class }$\mathcal{S}{}{}$\textit{ of
smooth projective varieties\footnote{By a variety, I mean a geometrically
reduced scheme of finite type over the ground field.} over }$k$,\textit{
closed under the formation of products, disjoint sums, and passage to a
connected component, and containing the abelian varieties, projective spaces,
and varieties of dimension zero. Except in the last section, we assume that
the Tate conjecture holds for the varieties in }$\mathcal{S}{}$\textit{.}

\section{The cohomology of groups of multiplicative type}

Let $M$ be a finitely generated $\mathbb{Z}{}$-module with a continuous action
of $\Gamma=\Gal(\mathbb{Q}{}^{\mathrm{al}}/\mathbb{Q}{})$ (discrete topology
on $M$), and let $T=D(M)$ be the corresponding algebraic group of
multiplicative type over $\mathbb{Q}{}$. Thus
\[
X^{\ast}(T)\overset{\text{{\tiny def}}}{=}\Hom(T_{\mathbb{Q}{}^{\mathrm{al}}%
},\mathbb{G}_{m})=M.
\]
For $m\in M$, let $\mathbb{Q}{}[m]$ be the fixed field of $\Gamma_{m}%
\overset{\text{{\tiny def}}}{=}\{\sigma\in\Gal(\mathbb{Q}^{\mathrm{al}%
}/\mathbb{Q}{})\mid\sigma m=m\}$, so that
\[
\Sigma_{m}\overset{\text{{\tiny def}}}{=}\Hom(\mathbb{Q}{}[m],\mathbb{Q}%
{}^{\mathrm{al}})\simeq\Gamma/\Gamma_{m}.
\]
Let $(\mathbb{G}_{m})_{\mathbb{Q}{}[m]/\mathbb{Q}{}}$ be the torus over
$\mathbb{Q}{}$ obtained from $\mathbb{G}_{m}$ by (Weil) restriction of scalars
from $\mathbb{Q}{}[m]$ to $\mathbb{Q}{}$, so that
\[
X^{\ast}((\mathbb{G}_{m})_{\mathbb{Q}{}[m]/\mathbb{Q}{}})\simeq\mathbb{Z}%
{}[\Sigma_{m}]
\]
(free $\mathbb{Z}{}$-module on $\Sigma_{m}$ with $\tau\in\Gamma$ acting by
$\tau(\sum n_{\sigma}\sigma)=\sum n_{\sigma}\tau\circ\rho$). The map
\begin{equation}
\tstyle\sum_{\sigma}n_{\sigma}\sigma\mapsto\tstyle\sum_{\sigma}n_{\sigma}%
\cdot\sigma m\colon\mathbb{Z}{}[\Sigma_{m}]\rightarrow M, \label{e12}%
\end{equation}
defines a homomorphism $T\rightarrow\mathbb{(}\mathbb{G}_{m})_{\mathbb{Q}%
{}[m]/\mathbb{Q}{}}$ and hence a homomorphism%
\begin{equation}
\alpha_{m}\colon H^{2}(\mathbb{Q}{},T)\rightarrow H^{2}(\mathbb{Q}%
{},\mathbb{(}\mathbb{G}_{m})_{\mathbb{Q}{}[m]/\mathbb{Q}{}})\simeq
\Br(\mathbb{Q}{}[m]). \label{e16}%
\end{equation}

\begin{proposition}
\label{1}Let $c\in H^{2}(\mathbb{Q}{},T)$. If $\alpha_{m}(c)=0$ for all $m\in
M$, then $c$ lies in the kernel of%
\[
H^{2}(\mathbb{Q}{},T)\rightarrow H^{2}(\mathbb{Q}{}_{l},T)
\]
for every prime $l$ of $\mathbb{Q}{}$ (including $l=\infty$).
\end{proposition}

\begin{proof}
Let $c$ be an element of $H^{2}(\mathbb{Q}{},T)$ such that $\alpha_{m}(c)=0$
for all $m\in M$, and fix a finite prime $l$ of $\mathbb{Q}{}$. To show that
$c$ maps to zero in $H^{2}(\mathbb{Q}{}_{l},T)$, it suffices to show that the
family of homomorphisms
\[
\alpha_{m,l}\colon H^{2}(\mathbb{Q}_{l},T)\rightarrow H^{2}(\mathbb{Q}%
_{l},(\mathbb{G}_{m})_{\mathbb{Q}[m]/\mathbb{Q}}),\quad m\in M,
\]
is injective. Choose an extension of $l$ to $\mathbb{Q}{}^{\mathrm{al}}$, and
let $\Gamma(l)\subset\Gamma$ be the corresponding decomposition group. A
standard duality theorem (\cite{milne1986ar}, I 2.4) shows that the
$\alpha_{m,l}$ is obtained from the homomorphism
\begin{equation}
\mathbb{Z}{}[\Sigma_{m}]^{\Gamma(l)}\rightarrow M^{\Gamma(l)} \label{e13}%
\end{equation}
by applying the functor $\Hom(\cdot,\mathbb{Q}/\mathbb{Z})$. Thus it suffices
to prove that the family of homomorphisms (\ref{e13}), indexed by $m\in M$, is
surjective. Let $m\in M^{\Gamma(l)}$. Because the group $\Gamma(l)$ fixes $m$,
it is contained in $\Gamma_{m}$, and so it fixes the inclusion $\sigma
_{0}\colon\mathbb{Q}[m]\hookrightarrow\mathbb{Q}^{\text{al}}$. Thus
$\sigma_{0}$ is an element of $\mathbb{Z}{}[\Sigma_{m}]^{\Gamma(l)}$, and it
maps to $m$.

The proof with $l=\infty$ is similar (apply \cite{milne1986ar}, I 2.13b).
\end{proof}

\begin{nt}
The proposition is abstracted from \cite{milne1994} (proof of Theorem 3.13).
\end{nt}

\section{Review of the category of motives over $\mathbb{F}$}

Let $q=p^{n}$. Recall that a \emph{Weil }$q$\emph{-number of weight }$m$ is an
algebraic number $\pi$ such that

\begin{itemize}
\item $|\pi^{\prime}|=q^{m/2}$ for all conjugates $\pi^{\prime}$ of $\pi$ in
$\mathbb{C}{}$ and

\item $q^{N}\pi$ is an algebraic integer for some $N$.
\end{itemize}

\noindent The first condition implies that $\pi\mapsto q^{m}/\pi$ defines an
automorphism $\iota^{\prime}$ of $\mathbb{Q}{}[\pi]$ such that $\sigma
\circ\iota^{\prime}=\iota\circ\sigma$ for all $\sigma\colon\mathbb{Q}{}%
[\pi]\rightarrow\mathbb{C}{}$. Therefore, $\mathbb{Q}{}[\pi]$ is totally real
or CM. Note that, because $q^{N}\pi$ is an algebraic integer and $(q^{N}%
\cdot\pi)(q^{N}\cdot\iota^{\prime}\pi)=q^{2N+m}$, the ideal $(\pi)$ is
divisible only by $p$-adic primes.

Fix a CM-subfield $K$ of $\mathbb{Q}^{\mathrm{al}}$, finite and galois over
$\mathbb{Q}{}$, and let $W_{0}^{K}(q)$ denote the set of Weil $q$-numbers of
weight $0$ in $K$ such that%
\[
n_{w}(\pi)\overset{\text{{\tiny def}}}{=}\frac{\ord_{w}(\pi)}{\ord_{w}%
(q)}[K_{w}\colon\mathbb{Q}{}_{p}]
\]
lies in $\mathbb{Z}{}$ for all $p$-adic primes $w$ of $K$. Note that the
torsion subgroup of $W_{0}^{K}(q)$ is $\mu(K)$, the group of roots of $1$ in
$K$. Let $X$ and $Y$ be the sets of $p$-adic primes of $K$ and of its largest
real subfield $F$. Write $f_{K}$ for the common inertia degree\footnote{The
inertia degree of a prime $\mathfrak{p}{}$ of $K$ is the degree
$f(\mathfrak{p}{}/p)=[\mathcal{O}{}_{K}/\mathfrak{p}{}\colon\mathbb{F}{}_{p}]$
of the field extension $\mathcal{O}{}_{K}/\mathfrak{p}{}\supset\mathbb{F}%
{}_{p}$.} of the $p$-adic prime ideals of $K$ and $h_{K}$ for their common
order in the class group of $K$.

\begin{proposition}
\label{4}For any $n$ divisible by $f_{K}h_{K}$, the sequence%
\begin{equation}
\minCDarrowwidth20pt\begin{CD} 0@>>>W_{0}^{K}(p^{n})/\mu (K)@>\pi\mapsto{\sum _{w|p}n_{w}(\pi)w}>>\mathbb{Z}^{X} @>{\sum a_{w}w\mapsto}{\sum a_{w}w|F}>>\mathbb{Z}{}^{Y}@>>>0 \end{CD} \label{e11}%
\end{equation}
is exact.
\end{proposition}

\begin{proof}
Everything is obvious except that every element in the kernel of the second
map is in the image of the first.

Let $\Gamma=\Hom(K,\mathbb{\mathbb{Q}{}}^{\mathrm{al}})=\Gal(K/\mathbb{Q}%
{}^{\mathrm{al}})$. The group $I_{0}(K)$ of infinity types of weight $0$ on
$K$ is the subgroup of $\mathbb{Z}{}[\Gamma]$ consisting of the sums $\sum
n_{\sigma}\sigma$ such that $n_{\sigma}+n_{\iota\sigma}=0$ for all $\sigma$.
Fix a $p$-adic prime $w_{0}$ of $\mathbb{Q}{}^{\mathrm{al}}$. As $\Gamma$ acts
transitively on $X$, the sequence%
\[
\begin{CD}
I_{0}(K)@>{\sum n_{\sigma}\sigma \mapsto \sum n_{\sigma}\sigma w_0}>>\mathbb{Z}%
{}^{X}@>{\sum a_{w}w\mapsto\sum a_{w}w|F}>>\mathbb{Z}^{Y}%
@>>>0
\end{CD}
\]
is exact. Because $n/f_{K}$ is divisible by $h_{K}$, there exists an element
$\varpi$ of $\mathcal{O}{}_{K}$ such that $\mathfrak{p}{}_{w_{0}}^{n/f_{K}%
}=(\varpi)$, i.e., such that for $w$ a finite prime of $K$,
\[
\ord_{w}(\varpi)=%
\begin{cases}
n/f_{K} & \text{if $w=w_{0}$}\\
0 & \text{otherwise.}%
\end{cases}
\]
For $\chi=\sum n_{\sigma}\sigma\in I_{0}(K)$ and $a\in K^{\times}$, let
\[
\chi(a)=\tsize\prod_{\sigma\in\Gamma}(\sigma a)^{n_{\sigma}}.
\]
Then $\chi(a)=1$ for $a\in F^{\times}$. As the group of units in $F$ has
finite index in the group of units in $K$, this shows that $\pi\overset
{\text{{\tiny def}}}{=}\chi(\varpi)$ is independent of the choice of $\varpi$
up to an element of $\mu(K)$. It lies in $W_{0}^{K}(p^{n})$, and the diagram%
\[
\xymatrixrowsep{0.5pc}\xymatrix{
I_{0}(K)\ar[rd]\ar[dd]^{\chi}\\
&\mathbb{Z}{}^{X}\\
W_{0}^{K}(p^{N})\ar[ru]}
\]
commutes, which completes the proof.
\end{proof}

If $n|n^{\prime}$, then $\pi\mapsto\pi^{n^{\prime}/n}$ is a homomorphism
$W^{K}(p^{n})\rightarrow W^{K}(p^{n^{\prime}})$; we define $W^{K}(p^{\infty
})=\varinjlim W^{K}(p^{n})$. Similarly, $W_{0}^{K}(p^{\infty})=\varinjlim
W_{0}^{K}(p^{n})$. Thus, an element of $W_{0}^{K}(p^{\infty})$ is represented
by a pair $(\pi,n)$ with $n\in\mathbb{N}{}^{\times}$ and $\pi\in W_{0}%
^{K}(p^{n})$.

\begin{corollary}
\label{5}When $n$ divides $n^{\prime}$ and both are divisible by $f_{K}h_{K}$,%
\[
W_{0}^{K}(p^{n})/\mu(K)\simeq W_{0}^{K}(p^{n^{\prime}})/\mu(K)\simeq W_{0}%
^{K}(p^{\infty}).
\]
The sequence%
\[
0\rightarrow W^{K}(p^{\infty})\rightarrow\mathbb{Z}{}^{X}\rightarrow
\mathbb{Z}{}^{Y}\rightarrow0
\]
is exact.
\end{corollary}

\begin{proof}
The diagram%
\[
\begin{CD}
W_{0}^{K}(p^{n}) @>>> \phantom{x}\mathbb{Z}^{X}\\
@V{\pi\mapsto\pi^{n/n^{\prime}}}VV@|\\
W_{0}^{K}(p^{n^{\prime}})@>>>\phantom{x}\mathbb{Z}^{X}%
\end{CD}
\]
commutes, and so this follows from (\ref{e11}).
\end{proof}

Let $P^{K}(p^{n})$ be the algebraic group of multiplicative type over
$\mathbb{Q}{}$ with character group $X^{\ast}(P^{K}(p^{n}))=W^{K}(p^{n})$.

\begin{corollary}
\label{6}Let $m=|\mu(K)|$. For any $n$ divisible by $f_{K}h_{K}$, there are
exact sequences%
\[
\minCDarrowwidth15pt \begin{CD}
0 @>>> (\mathbb{G}_{m})_{F/\mathbb{Q}} @>>> (\mathbb{G}_{m})_{K/\mathbb{Q}}@>>>P^{K}(p^{\infty})@>>>0\\
@.       @.                                      @.                                @|\\
@.@.0 @>>> P^{K}(p^{\infty}) @>>> P^{K}(p^{n}) @>>> \mathbb{Z}/m\mathbb{Z} @>>> 0
\end{CD}
\]

\end{corollary}

\begin{proof}
Obvious from (\ref{4}) and (\ref{5}).
\end{proof}

\begin{proposition}
\label{7}(a) The family of maps $H^{1}(\mathbb{Q}{},P_{0}^{K}(p^{\infty
}))\rightarrow H^{1}(\mathbb{Q}{}_{l},P_{0}^{K}(p^{\infty})$, with $l$ running
over the primes of $\mathbb{Q}{}$, is injective.

(b) $H^{2}(\mathbb{Q}{},P_{0}^{K}(p^{\infty}))\simeq\bigoplus\nolimits_{l}%
H^{2}(\mathbb{Q}{}_{l},P_{0}^{K}(p^{\infty}))$ (sum over all primes of
$\mathbb{Q}{}$).
\end{proposition}

\begin{proof}
This follows from the cohomology sequence of the upper exact sequence in
(\ref{6}) and class field theory (\cite{milne1994}, 3.11).
\end{proof}

\begin{proposition}
\label{9}The family of maps $H^{2}(\mathbb{Q}{},P_{0}^{K}(p^{\infty
}))\rightarrow\Br(\mathbb{Q}{}[\pi])$, $\pi\in W_{0}^{K}(p^{\infty})$, (see
(\ref{e16})) is injective.
\end{proposition}

\begin{proof}
Apply Proposition \ref{1} and (b) of Proposition \ref{7}.
\end{proof}

Let $X$ be a simple object in $\Mot_{0}^{K}(\mathbb{F}{})$, and let $\pi_{X}$
be its Frobenius endomorphism.

\begin{proposition}
\label{8}The map $H^{2}(\alpha_{\pi_{X}})\colon H^{2}(\mathbb{Q}{},P_{0}%
^{K}(p^{\infty}))\rightarrow\Br(\mathbb{Q}{}[\pi_{X}])$ sends the class of
$\Mot_{0}^{K}(\mathbb{F}{})$ in $H^{2}(\mathbb{Q}{},P_{0}^{K}(p^{\infty}))$ to
the class of $\End(X)$ in $\Br(\mathbb{Q}{}[\pi_{X}])$.
\end{proposition}

\begin{proof}
This can be proved by the same argument as in \cite{saavedra1972}, VI 3.5.3.
\end{proof}

\begin{remark}
\label{9a}Let $W_{0}(p^{\infty})=\varinjlim\nolimits_{n}W_{0}(p^{n})$ and
$P_{0}(p^{\infty})=\varprojlim_{K}P_{0}^{K}(p^{\infty})$; thus $P_{0}%
(p^{\infty})$ is the pro-torus with character group $W_{0}(p^{\infty})$.
Proposition \ref{9} shows that the family of maps $H^{2}(\mathbb{Q}{}%
,P_{0}(p^{\infty}))\rightarrow\Br(\mathbb{Q}{}[\pi])$, $\pi\in W_{0}%
(p^{\infty})$, has kernel $\varprojlim_{K}^{1}H^{1}(\mathbb{Q}{},P_{0}%
^{K}(p^{\infty}))$, which is zero (\cite{milne2003}, 3.8).
\end{remark}

\begin{nt}
This section reviews results from \cite{langlandsR1987, wei1993, milne1994,
milne2003}.
\end{nt}

\section{The category of motives over $\mathbb{F}{}_{p}$}

As before, fix a (large) CM-subfield $K$ of $\mathbb{Q}^{\mathrm{al}}$, finite
and galois over $\mathbb{Q}{}$. Let $W_{0}^{K}(p,n)$ be the group of Weil
$p$-numbers $\pi$ in $\mathbb{Q}{}^{\mathrm{al}}$ of weight $0$ such that
$\pi^{n}\in W_{0}^{K}(p^{n})$. Note that, $W_{0}^{K}(p,1)=W_{0}^{K}(p)$, but
otherwise the elements of $W_{0}^{K}(p,n)$ need not lie in $K$.

\begin{lemma}
\label{10}For any $n$ divisible by $f_{K}h_{K}$, there is an exact sequence%
\[
\begin{CD}
0@>>>\mu_{mn}@>>> W_{0}^{K}(p,n)@>{\pi\mapsto\lbrack\pi^{n},n]}>>W_{0}^{K}(p^{\infty})@>>>0.
\end{CD}
\]
where $m=|\mu(K)|$ and $\mu_{mn}=\mu_{mn}(\mathbb{Q}{}^{\mathrm{al}})$.
\end{lemma}

\begin{proof}
According to (\ref{5}), an element of $W_{0}^{K}(p^{\infty})$ is represented
by a $\pi\in W_{0}^{K}(p^{n})$. Now any $n$th root $\pi^{\frac{1}{n}}$ of
$\pi$ in $\mathbb{Q}{}^{\mathrm{al}}$ lies in $W_{0}^{K}(p,n)$ and maps to
$\pi$.

If $\pi\in W_{0}^{K}(p,n)$ is such that $\pi^{n}$ represents $1$ in $W_{0}%
^{K}(p^{\infty})$, then $\pi^{n}\in W_{0}^{K}(p^{n})_{\text{tors}}=\mu(K)$.
Therefore $(\pi^{n})^{m}=1$. Conversely, if $\pi^{mn}=1$, then $\pi\in
W_{0}(p,n)$ and maps to $1$ in $W_{0}^{K}(p^{\infty})$.
\end{proof}

Let $P_{0}^{K}(p,n)$ be the group of multiplicative type over $\mathbb{Q}{}$
such that $X^{\ast}(P_{0}^{K}(p,n))=W_{0}^{K}(p,n)$.

\begin{proposition}
\label{11}For any $n$ divisible by $f_{K}h_{K}$, there is an exact sequence%
\begin{equation}
0\rightarrow P_{0}^{K}(p^{\infty})\rightarrow P_{0}^{K}(p,n)\rightarrow
\mathbb{Z}/nm\mathbb{Z}{}\rightarrow0. \label{e4}%
\end{equation}

\end{proposition}

\begin{proof}
Immediate from the lemma.
\end{proof}

Recall that the isomorphism classes of simple objects in $\Mot(\mathbb{F}%
{}_{p})$ are classified by the conjugacy classes of elements of $W(p)$ (Weil
$p$-numbers in $\mathbb{Q}{}^{\mathrm{al}}$) (see, for example,
\cite{milne1994}, 2.6). Let $\Mot_{0}^{K}(\mathbb{F}{}_{p},n)$ be the category
of motives over $\mathbb{F}{}_{p}$ whose Weil $p$-numbers lie in $W_{0}%
^{K}(p,n)$.

\begin{proposition}
\label{11a}The category $\Mot_{0}^{K}(\mathbb{F}{}_{p},n)$ has a $\mathbb{Q}%
{}_{l}$-valued fibre functor for all primes $l$ of $\mathbb{Q}{}$ (including
$p$ and $\infty$).
\end{proposition}

\begin{proof}
As we noted in the introduction, the endomorphism algebras of simple objects
in $\Mot_{0}^{K}(\mathbb{F}{}_{p},n)$ are commutative, and so this follows
from Propositions \ref{8} and \ref{1}.
\end{proof}

\section{Cyclic algebras}

Let $F$ be field.

\begin{definition}
\label{22}A \emph{cyclic semifield over }$F$ is an \'{e}tale $F$-algebra $E$
together with an action of a cyclic group $C$ such that $C$ acts simply
transitively on $\Hom_{F\text{-algebra}}(E,F^{\mathrm{al}})$. In other words,
it is a galois $F$-algebra with cyclic galois group (in the sense of Grothendieck).
\end{definition}

\begin{proposition}
\label{23}Let $E$ be a cyclic field extension of $F$ with generating
automorphism $\sigma_{0}$. Then $(E^{m},\sigma)$, with $\sigma(a_{1}%
,\ldots,a_{m})=(\sigma_{0}a_{m},a_{1},\ldots,a_{m-1})$, is a cyclic semifield
over $F$, and every cyclic semifield over $F$ is isomorphic to one of this form.
\end{proposition}

\begin{proof}
Routine application of galois theory (in the sense of Grothendieck).
\end{proof}

We denote $(E^{m},\sigma)$ by $(E,\sigma_{0})^{m}$.

\begin{example}
\label{23a}Let $(E,\sigma)$ be a cyclic field over $F$. Let $F^{\prime}$ be a
field containing $F$, and let $EF^{\prime}$ be the composite of $E$ and
$F^{\prime}$ in some common larger field. Let $m$ be the least positive
integer such that $\sigma^{m}$ fixes $E\cap F^{\prime}$. Then%
\[
\langle\sigma^{m}\rangle=\Gal(E/E\cap F^{\prime})\simeq\Gal(EF^{\prime
}/F^{\prime}),
\]
and so $(EF^{\prime},\sigma^{m})$ is a cyclic field over $F$. Clearly,
$(E\otimes_{F}F^{\prime},\sigma)\approx(EF^{\prime},\sigma^{m})^{m}$.
\end{example}

Let $(E,\sigma)$ be a cyclic semifield over $F$. For any element $a\in
F^{\times}$, define%
\[
B(E,\sigma,a)=E\cdot1+E\cdot x+\cdots+Ex^{n-1}%
\]
with the multiplication determined by%
\[
x^{n}=a,\quad x\cdot e=\sigma(e)\cdot x\text{ for }b\in E.
\]
Then $B(E,\sigma,a)$ is a central simple $F$-algebra (\cite{albert1939}, VII).
Algebras of this form are called \emph{cyclic}.\footnote{Classically, they
were called \textquotedblleft generalized cyclic algebras\textquotedblright,
and \textquotedblleft cyclic algebra\textquotedblright\ was reserved for those
with $E$ is a field.} Because $E$ is a maximal \'{e}tale subalgebra of $B$, it
splits $B$, and so $B(E,\sigma,a)$ represents a class in $\Br(E/F)$.

\begin{proposition}
\label{24}With the notations of (\ref{23}),
\[
B(E^{m},\sigma,a)\approx B(E,\sigma_{0},a)\otimes_{F}M_{m}(F).
\]

\end{proposition}

\begin{proof}
See \cite{albert1939}, VII 1, Theorem 1.
\end{proof}

\begin{corollary}
\label{24a}Let $(E,\sigma)$ be a cyclic field over $F$. Let $F^{\prime}$ be a
field containing $F$, and let $a\in F^{\times}$. With the notations of
(\ref{23a}),
\[
\lbrack B(E\otimes_{F}F^{\prime},\sigma,a)]=[B(E,\sigma,a)\otimes_{F}%
F^{\prime}]=[B(EF^{\prime},\sigma^{m},a)]
\]
(equality of classes in $\Br(F^{\prime})$). More generally, when $a\in
F^{\prime\times}$, one still has%
\[
\lbrack B(E\otimes_{F}F^{\prime},\sigma,a)]=[B(EF^{\prime},\sigma^{m},a)]
\]

\end{corollary}

\begin{proof}
Apply the proposition to $(E\otimes_{F}F^{\prime},\sigma)\approx(EF^{\prime
},\sigma^{m})^{m}$. See also \cite{reiner2003}, 30.8, for the case where $a\in
F^{\times}$.
\end{proof}

\begin{plain}
\label{12a}For a fixed $(E/F,\sigma)$, the map
\[
a\mapsto\lbrack B(E,\sigma,a)]\colon F^{\times}\rightarrow\Br(E/F)
\]
has the following cohomological description. The choice of the generator
$\sigma$ for the galois group of $E/F$ determines an isomorphism of the Tate
cohomology groups $H^{0}(E/F,E^{\times})\rightarrow H^{2}(E/F,E^{\times})$,
i.e., an isomorphism%
\begin{equation}
F^{\times}/\Nm E^{\times}\rightarrow\Br(E/F) \label{e5}%
\end{equation}
(periodicity of the cohomology of cyclic groups; see, for example,
\cite{milneCFT}, II 2.11). This isomorphism maps $a\in F^{\times}$ to the
class of $B(E,\sigma,a)$. When $F$ is a local or global field, it is known
that every element of $\Br(F)$ is split by a cyclic extension, and so is
represented by a cyclic algebra.
\end{plain}

\begin{example}
\label{12}Let $F$ be a finite extension of $\mathbb{Q}{}_{p}$ and let $E$ be
an unramified field extension of $F$ of degree $n$. Choose $\sigma$ to be the
Frobenius element. For any $a\in F^{\times}$, $B(E,\sigma,a)$ has invariant
$\ord_{K}(a)/n$ (cf. \cite{milneCFT}, IV 4.2). Here $\ord_{K}$ is normalized
to map onto $\mathbb{Z}{}$.
\end{example}

We now fix $(F,a)$ and give a cohomological description of
\begin{equation}
(E,\sigma)\mapsto\lbrack B(E,\sigma,a)]\colon H^{1}(F,\mathbb{Z}{}%
/n\mathbb{Z}{})\rightarrow\Br(F). \label{e1}%
\end{equation}
Let $F[x]=F[X]/(X^{n}-a)$. The inclusion $F^{\times}\hookrightarrow
F[x]^{\times}$ defines a homomorphism $\mathbb{G}_{m}\rightarrow
(\mathbb{G}_{m})_{F[x]/F}$, and we let $T$ be the cokernel. The class of $x$
in $F[x]^{\times}/F^{\times}\subset T(\mathbb{Q}{})$ has order dividing $n$,
and the diagram%
\[
\begin{CD}
1 @>>> \mathbb{G}_{m} @>>> (\mathbb{G}_{m})_{F[x]/F} @>>> T
@>>> 0\\
@.@.@.@AA1\mapsto[x]A \\
@.@.@.\mathbb{Z}/n\mathbb{Z}
\end{CD}
\]
of groups of multiplicative type gives rise to a diagram of cohomology groups%
\[
\begin{CD}
H^{1}(F,T) @>>> \Br(F)@>>>\Br(F[x])\\
@AAA\\
H^{1}(F,\mathbb{Z}{}/n\mathbb{Z})
\end{CD}
\]

\begin{lemma}
\label{13}The composite of the maps%
\[
H^{1}(F,\mathbb{Z}{}/n\mathbb{Z}{})\rightarrow H^{1}(F,T)\rightarrow\Br(F)
\]
is the map (\ref{e1}).
\end{lemma}

\begin{proof}
Omitted (for the moment).
\end{proof}

\begin{proposition}
\label{14}Let $L$ be an unramified cyclic field extension of $\mathbb{Q}{}%
_{p}$ of degree $n$, and let $\sigma$ be the Frobenius automorphism of $L$
over $\mathbb{Q}{}_{p}$. For any finite extension $F$ of $\mathbb{Q}{}_{p}$
and $a\in F$,
\begin{equation}
\inv_{v}B(L\otimes_{\mathbb{Q}{}_{p}}F,\sigma,a)=\frac{\ord_{F}(a)}%
{n\cdot\ord_{F}(p)}\cdot\lbrack F\colon\mathbb{Q}{}_{p}]. \label{e14}%
\end{equation}

\end{proposition}

\begin{proof}
Let $e$ and $f$ be the ramification and inertia indices of $p$ in $F$. The
composite $LF$ of $L$ and $F$ in some common larger field is an unramified
extension of $F$ of degree $n/f$ with Frobenius element $\sigma^{f}$. Thus
(see \ref{24a}, \ref{12}),
\begin{align*}
\inv_{F}B(L\otimes_{\mathbb{Q}{}_{p}}F,\sigma,a))  &  =\inv_{F}(B(LF,\sigma
^{f},a)\\
&  =\frac{\ord_{F}(a)}{n/f}.
\end{align*}
Since $\ord_{F}(p)=e$, this gives (\ref{e14}).
\end{proof}

\section{The $\mathbb{Q}{}$-valued fibre functors on $\Mot_{0}^{K}%
(\mathbb{F}{}_{p},n)$}

As before, $K$ is a CM subfield of $\mathbb{Q}{}^{\mathrm{al}}$, finite and
galois over $\mathbb{Q}{}$, and $n$ is an integer divisible by $f_{K}h_{K}$.

\begin{lemma}
\label{16a}Let $\pi\in W_{0}^{K}(p,n)$, and let $\bar{\pi}$ be the image of
$\pi$ in $W_{0}(p^{\infty})$. Let $X(\bar{\pi})$ be the simple motive over
$\mathbb{F}{}$ corresponding to $\bar{\pi}$. Then the centre of $\End(X(\bar
{\pi}))$ is $\mathbb{Q}{}[\pi^{mn}]$.
\end{lemma}

\begin{proof}
Recall that $\bar{\pi}$ is represented by $\pi^{n}\in W_{0}^{K}(p^{n})$. The
centre of $\End(X(\bar{\pi}))$ is $\mathbb{Q}{}[\bar{\pi}]$ (notations as in
\S 1, i.e., $\mathbb{Q}{}[\bar{\pi}]$ is the fixed field of the subgroup of
$\Gal(\mathbb{Q}{}^{\mathrm{al}}/\mathbb{Q}{})$ fixing $\bar{\pi}$). An
element $\sigma$ of $\Gal(\mathbb{Q}{}^{\mathrm{al}}/\mathbb{Q}{})$ fixes
$\bar{\pi}$ if and only if it fixes its image in $\mathbb{Z}{}^{X}$ (notation
as in \S 2), but this equals the image $\pi^{n}$ in $\mathbb{Z}{}^{X}$, which
is fixed by $\sigma$ if and only if $\sigma$ fixes $(\pi^{n})^{m}$.
\end{proof}

\begin{theorem}
\label{17}There exists a $\mathbb{Q}{}$-valued fibre functor $\omega$ on
$\Mot_{0}^{K}(p,n)$ if and only if there exists a cyclic field extension $L$
of $\mathbb{Q}{}$ of degree $mn$ such that

\begin{enumerate}
\item $(p)$ remains prime in $L$;

\item $\pi^{mn}$ is a local norm at every prime $v$ of $\mathbb{Q}{}[\pi
^{mn}]$ that ramifies in $\mathbb{Q}{}[\pi^{mn}]\otimes_{\mathbb{Q}{}}L$.
\end{enumerate}
\end{theorem}

\begin{proof}
Consider the diagram arising from (\ref{e4}) and (\ref{e1})
\[
\xymatrix{
H^{1}(\mathbb{Q},\mathbb{Z}/mn\mathbb{Z})\ar[r]^{\alpha}\ar[rd]_{\gamma}
&H^{2}(\mathbb{Q},P_{0}^{K}(p^{\infty})\ar[r]^{\beta}\ar[d]^{\textrm{injective}}
&H^{2}(\mathbb{Q},P_{0}^{K}(p,n))\\
&\displaystyle\prod_{\pi\in W_0^K(p,n)}\Br(\mathbb{Q}{}[\pi^{mn}])
}
\]
The map $\beta$ sends the cohomology class of $\Mot_{0}^{K}(\mathbb{F}{}_{p})$
to that of $\Mot_{0}^{K}(\mathbb{F}{}_{p},n)$. Thus, $\Mot_{0}^{K}%
(\mathbb{F}{}_{p},n)$ is neutral if and only if the cohomology class of
$\Mot_{0}^{K}(\mathbb{F}{}_{p})$ is in the image of $\alpha$. Since $\gamma$
sends an element $(L,\sigma)$ of $H^{1}(\mathbb{Q}{},\mathbb{Z}{}%
/mn\mathbb{Z}{})$ to the class of $B(L\otimes_{\mathbb{Q}{}}\mathbb{Q}{}%
[\pi^{mn}],\sigma,\pi^{mn})$ in $\Br(\mathbb{Q}{}[\pi^{mn}])$ (cf. \ref{13}),
we see that $\Mot_{0}^{K}(p,n)$ is neutral if and only if there exists a
cyclic field extension $(L,\sigma)$ of degree dividing $mn$ such that, for all
$\pi\in W_{0}^{K}(p,n)$ and all primes $v$ of $\mathbb{Q}{}[\pi^{mn}]$,%
\begin{equation}
\inv_{v}(B(L\otimes_{\mathbb{Q}{}}\mathbb{Q}{}[\pi^{mn}],\sigma,\pi
^{mn}))=\left\{
\begin{array}
[c]{ll}%
\dfrac{\ord_{v}(\pi)}{\ord_{v}(p^{mn})}\cdot\lbrack\mathbb{Q}{}[\pi^{mn}%
]_{v}:\mathbb{Q}_{p}] & \text{if }v|p\\
0 & \text{otherwise}%
\end{array}
\right.  \label{e15}%
\end{equation}

Let $L$ be a cyclic field extension of $\mathbb{Q}{}$ of degree $mn$
satisfying the conditions (a) and (b) and let $\sigma=(p,L/\mathbb{Q}{})$.
Condition (a) implies that (\ref{e15}) holds for the primes $v$ dividing $p$
(apply Proposition \ref{14} with $L\otimes_{\mathbb{Q}{}}\mathbb{Q}{}_{p}$ for
$L$, $\mathbb{Q}{}[\pi^{mn}]_{v}$ for $F$, $\pi^{mn}$ for $a$), and condition
(b) implies that (\ref{e15}) holds for the primes not dividing $p$ (see
(\ref{12a})).

Conversely, let $(L,\sigma)$ be a cyclic extension of $\mathbb{Q}{}$ of degree
dividing $mn$ satisfying (\ref{e15}). By considering the primes dividing $p$
and applying Proposition \ref{14}, one sees that $L$ has degree $mn$, that
$(p)$ remains primes in $L$, and $\sigma=(p,L/K)$. On the other hand, the
invariant at a prime not dividing $p$ vanishes automatically unless the prime
ramifies in $\mathbb{Q}{}[\pi^{mn}]$, in which case it vanishes if and only if
(b) holds (by \ref{12a}).
\end{proof}

\begin{theorem}
\label{18}Let $F$ be a the field generated over $\mathbb{Q}{}$ by the elements
of $W_{0}^{K}(p,n)$ --- it is a finite galois extension of $\mathbb{Q}{}$. The
generalized Riemann hypothesis implies that there exists a field $L$
satisfying the conditions (a) and (b) of (\ref{17}) provided $p$ is not an
$r$th power in $F\cdot\mathbb{Q}{}^{\mathrm{ab}}$ for any $r$ dividing $mn$.
\end{theorem}

\begin{proof}
Note that $F$ is generated over $\mathbb{Q}{}$ by any set of generators for
the abelian group $W_{0}^{K}(p,n)$, which can be chosen to be finite and
stable under the action of $\Gal(\mathbb{Q}{}^{\mathrm{al}}/\mathbb{Q}{})$,
which shows that $F$ is finite and galois over $\mathbb{Q}{}$. Note that
condition (b) is implied by the stronger condition: \bquote(b$^{\prime}$)
every prime $l\neq p$ ramifying in $L$ splits in $F.$\equote\noindent
\ (Because then $\pi\in\mathbb{Q}{}[\pi^{mn}]_{v}=\mathbb{Q}{}_{l}$, and so
$\pi^{mn}$ is an $mn$th power \textit{inside }$\mathbb{Q}{}[\pi^{mn}]_{v}$.) A
natural place to look for such an extension $L$ is inside $\mathbb{Q}{}%
[\zeta_{l}]$ for some prime $l$. Since only $l$ ramifies in $\mathbb{Q}%
{}[\zeta_{l}]$, it will contain an $L$ satisfying (a) and (b$^{\prime}$) if

\begin{enumerate}
\item[(c)] $(\mathbb{Z}{}/l\mathbb{Z}{})^{\times}$ has a quotient of order
$mn$ generated by the class of $p$,

\item[(d)] $l$ splits in $F$.
\end{enumerate}

\noindent We show in the next section that the generalized Riemann hypothesis
implies that, under our hypothesis on $p$, $r$, and $F$, there are always
infinitely many primes satisfying these conditions (c,d).
\end{proof}

\begin{remark}
\label{19}We can make the relation between the $\mathbb{Q}{}$-valued fibre
functors on $\Mot_{0}^{K}(\mathbb{F}{}_{p},n)$ and the cyclic field extensions
of $\mathbb{Q}{}$ more precise. The base change functor%
\[
\beta\colon\Mot_{0}^{K}(\mathbb{F}{}_{p},n)\rightarrow\Mot_{0}^{K}%
(\mathbb{F}{})
\]
realizes the second category as a normal quotient of the first category (in
the sense of \cite{milne2005}, \S 2). The objects of $\Mot_{0}^{K}%
(\mathbb{F}{}_{p},n)$ becoming trivial in $\Mot_{0}^{K}(\mathbb{F}{},n)$ are
exactly the Artin motives. Let $\omega^{\beta}$ be the fibre functor on
$\Art^{K}(\mathbb{F}{}_{p},n)$ defined by $\beta$ (ib. \S 2). Note that the
fundamental group of $\Art^{K}(\mathbb{F}{}_{p},n)$ is $\mathbb{Z}%
{}/mn\mathbb{Z}{}$, and that the motive $X_{mn}$ of $\mathbb{F}{}_{p^{mn}}$
lies in $\Art^{K}(\mathbb{F}{}_{p},n)$.

Now let $\omega$ be a $\mathbb{Q}{}$-valued fibre functor on $\Mot_{0}%
^{K}(\mathbb{F}{}_{p},n)$, and let $\wp=\underline{\Hom}^{\otimes}%
(\omega|,\omega^{q})$ where $\omega|$ is the restriction of $\omega$ to
$\mathrm{\Art}^{K}(\mathbb{F}{}_{p},n)$. Then $\wp$ is a $\mathbb{Z}%
{}/mn\mathbb{Z}{}$-torsor whose class in $H^{1}(\mathbb{Q}{},\mathbb{Z}%
{}/mn\mathbb{Z}{})$ maps to the class of $\Mot_{0}^{K}(\mathbb{F}{})$ in
$H^{2}(\mathbb{Q}{},P_{0}^{K}(p^{\infty}))$ (ib. 2.11). On the other hand, one
sees easily that the class of $\wp$ in $H^{1}(\mathbb{Q}{},\mathbb{Z}%
{}/mn\mathbb{Z}{})$ is represented by $L=\omega(X_{mn})$.
\end{remark}

We have seen that each $\mathbb{Q}{}$-valued fibre functor $\omega$ on
$\Mot_{0}^{K}(\mathbb{F}{}_{p})$ gives rise to a cyclic extension
$L=\omega(X_{mn})$ of $\mathbb{Q}{}$, and we have characterized the cyclic
extensions that arise in this way. To complete the classification, we have to
describe the set of fibre functors giving rise to the same field.

\begin{theorem}
\label{20}Let $\omega$ be a $\mathbb{Q}{}$-valued fibre functor on
$\Mot_{0}^{K}(p,n)$. The isomorphism classes of pairs consisting of a
$\mathbb{Q}{}$-valued fibre functor $\omega^{\prime}$ and an isomorphism
$\omega(X_{mn})\rightarrow\omega^{\prime}(X_{mn})$ are classified by
$\Br(E/F)$ where $E$ is the fixed field of the decomposition group of a
$p$-adic prime of $K$ and $F$ is its largest real subfield.
\end{theorem}

\begin{proof}
Let $\wp(\omega^{\prime})$ be the set of isomorphisms $\omega\rightarrow
\omega^{\prime}$ inducing the given isomorphism on $X_{mn}$. Then $\wp
(\omega^{\prime})$ is a torsor for $P_{0}^{K}(p^{\infty})$ (cf.
\cite{milne2004}, 1.6), and $\wp(\omega^{\prime})\approx\wp(\omega
^{\prime\prime})$ if and only if $\omega^{\prime}\approx\omega^{\prime\prime}%
$. Therefore, the pairs modulo isomorphism are classified by $H^{1}%
(\mathbb{Q}{},P_{0}^{K}(p^{\infty}))$, which equals $\Br(E/F)$
(\cite{milne1994}, 3.10).
\end{proof}

\section{The existence of the field $L$}

Let $a\neq\pm1$ be a square-free integer, let $k$ be a second integer, and let
$F$ be a finite galois extension of $\mathbb{Q}{}$. Consider the set $M$ of
prime numbers $p$ such that

\begin{itemize}
\item $p$ does not divide $a$,

\item $p$ splits in $F$,

\item the index in $(\mathbb{Z}{}/p\mathbb{Z}{})^{\times}$ of the subgroup of
generated by the class of $a$ divides $k{}.$
\end{itemize}

\noindent For each prime number $l$, let $q(l)$ be the smallest power of $l$
not dividing $k$, and let $L_{l}=\mathbb{Q}{}[\zeta_{q(l)},a^{1/q(l)}]$ be the
splitting field of $X^{q(l)}-a$ over $\mathbb{Q}{}$. If $p$ does not divide
$a$, then%
\[
p\text{ splits in }L_{l}\iff\left\{
\begin{array}
[c]{l}%
l|p-1\text{, and}\\
a\text{ is a }q(l)\text{th power modulo }p\text{.}%
\end{array}
\right.
\]
Therefore, a necessary condition for $M$ to be nonempty is that none of the
fields $L_{l}$ be contained in $F$.

\begin{theorem}
\label{20a}If the generalized Riemann hypothesis holds for each field $L_{l}$
and no $L_{l}$ is contained in $F$, then the set $M$ is infinite.
\end{theorem}

\begin{proof}
When $k=1$ and $F=\mathbb{Q}{}$, the statement becomes Artin's primitive root
conjecture: every square-free integer $a\neq\pm1$ is a primitive root for
infinitely many prime numbers $p$. That this follows from the generalized
Riemann hypothesis for the fields $L_{l}$ was proved by \nocite{hooley1967}%
Hooley (1967). The general case is proved in \cite{lenstra1977},
4.6.\footnote{Note that Lenstra frequently muddles his quantifiers. For
example, his condition \textquotedblleft$a_{n}\neq0$ for all $n$%
\textquotedblright\ should read \textquotedblleft no $a_{n}$ is
zero\textquotedblright.}
\end{proof}

\begin{lemma}
\label{20b}Let $a\neq\pm1$ be a square-free integer, and let $F$ be a finite
galois extension of $\mathbb{Q}{}$. Then there exists an integer $N$ such
that, if $a$ is an $m$th power in $F\cdot\mathbb{Q}{}^{\mathrm{ab}}$, then
$m|N$.
\end{lemma}

\begin{proof}
For odd primes $l$, the galois group of $X^{l}-a$ is never commutative, and so
$a$ is not an $l$th power in $\mathbb{Q}{}^{\mathrm{ab}}$. It follows that,
for any odd $m$, $X^{m}-a$ is irreducible over $\mathbb{Q}{}^{\mathrm{ab}}$
(e.g., \cite{lang2002}, VI Theorem 9.1, p297). Therefore, if $a$ is an $m$th
power in $\mathbb{Q}{}^{\mathrm{ab}}$, then $m|[F\mathbb{Q}{}^{\mathrm{ab}%
}\colon\mathbb{Q}^{\mathrm{ab}}]$.

The proof for even $m$ is similar.
\end{proof}

\begin{theorem}
\label{20c}Let $a\neq\pm1$ be a square-free integer, let $n$ be a positive
integer, and let $F$ be a finite galois extension of $\mathbb{Q}{}$. Let $M$
be the set of prime numbers $p$ such that

\begin{itemize}
\item $p$ does not divide $a$,

\item $p$ splits in $F$, and

\item $(\mathbb{Z}{}/p\mathbb{Z}{})^{\times}$ has a quotient of order $n$
generated by the class of $a$.
\end{itemize}

The set $M$ is empty if $a$ is an $m$th power in $F$ for some $m>1$ dividing
$n$, and it is infinite if $a$ is not an $m$th power in $F\cdot\mathbb{Q}%
{}^{\mathrm{ab}}$ for any $m$ dividing $n$.
\end{theorem}

\begin{proof}
Suppose $p\in M$. If $a$ is an $m$th power in $F$ for some $m$ dividing $n$,
then, because $p$ splits in $F$, $a$ is an $m$th power in $\mathbb{Q}{}_{p}$.
Therefore, it is an $m$th power in $(\mathbb{Z}{}/p\mathbb{Z}{})^{\times}$,
and in any cyclic quotient $C_{n}$ of $(\mathbb{Z}{}/p\mathbb{Z}{})^{\times}$.
Therefore, it can't generate $C_{n}$.

For the converse statement, the condition on $a$ implies that there exists a
$k$ relatively prime to $n$ such that $a$ is not a $q(l)$th power in
$F\cdot\mathbb{Q}{}^{\mathrm{ab}}$ for any prime $l$ (with $q(l)$ defined as
above). Then none of the fields $L_{l}$ is contained in $F[\zeta_{kn}]$, and
so there exist infinitely many primes $p$ such that

\begin{itemize}
\item $p$ does not divide $a$,

\item $p$ splits in $F[\zeta_{nk}]$,

\item the index in $(\mathbb{Z}{}/p\mathbb{Z}{})^{\times}$ of the subgroup of
generated by the class of $a$ divides $k$.
\end{itemize}

\noindent Because $p$ splits in $\mathbb{Q}{}[\zeta_{kn}]$, $kn$ divides
$p-1$, and so $(\mathbb{Z}{}/p\mathbb{Z}{})^{\times}$ has a quotient $C_{n}$
of order $n$. Because $k$ and $n$ are relatively prime, the image of $a$ in
$C_{n}$ generates it.
\end{proof}

\begin{corollary}
\label{20d}Theorem \ref{18} holds.
\end{corollary}

\begin{proof}
Apply the theorem with $(a,n)$ replaced by $(p,mn)$.
\end{proof}

\begin{question}
\label{20e}Does there exist a field $L$ satisfying conditions (a) and (b) of
Theorem \ref{17} for each pair $(K,n)$?
\end{question}

I don't see how to remove the proviso in Theorem \ref{18} much less the appeal
to the generalized Riemann hypothesis. In fact, I suspect that the answer to
the question is no. Here are two comments:

\begin{itemize}
\item Let $L=\mathbb{Q}{}[\zeta_{l}]$, and let $a\in\mathbb{Z}{}$ be
relatively prime to $l$. When is $a$ local norm at $l$? As $l$ is totally
ramified in $L$, the local Galois group is $(\mathbb{Z}{}/l\mathbb{Z}%
{})^{\times}$, and so this is true if and only if $a\equiv1$ modulo $l$.
Similarly, $a$ is a local norm from the subextension of $\mathbb{Q}{}%
[\zeta_{l}]$ of degree $m$ if and only if $a$ is an $m$th power in
$(\mathbb{Z}{}/l\mathbb{Z}{})^{\times}$.

\item See \cite{wei1993} for a description of the subfields of a CM-field
generated by Weil numbers.
\end{itemize}

\section{Fibre functors on $\Mot^{K}(\mathbb{F}{}_{p},n)$}

\begin{proposition}
\label{21}If there exists a $\mathbb{Q}{}$-valued fibre functor on
$\Mot_{0}^{K}(\mathbb{F}{}_{p},n)$, then
\end{proposition}

\begin{enumerate}
\item there exists a $\mathbb{Q}{}$-valued fibre functor on $\Mot_{\text{even}%
}^{K}(\mathbb{F}{}_{p},n)$, and

\item for a number field $L$, there exists an $L$-valued fibre functor on
$\Mot^{K}(\mathbb{F}{}_{p},n)$ if and only if the local degrees of the real
and $p$-adic primes of $L$ are even.
\end{enumerate}

\begin{proof}
Omitted (for the present).
\end{proof}

\section{Explicit description of the categories of motives}

In this section, we assume there exists an $L$ as in Theorem and give explicit
descriptions of various categories of motives.

\subsection{The category $\Mot_{0}^{K}(\mathbb{F}{}_{p},n)$}

The choice of a fibre functor $\omega$ on defines an equivalence
$X\mapsto\omega(X)$ from $\Mot_{0}^{K}(\mathbb{F}{}_{p},n)$ to the tannakian
category whose objects are the pairs $(V,F)$ with $V$ a finite-dimensional
vector space over $\mathbb{Q}{}$ and $F$ a semisimple endomorphism of $V$
whose eigenvalues lie in $W_{0}^{K}(p,n)$.

\subsection{The category $\Mot_{0}^{K}(\mathbb{F}{})$}

The realization of $\Mot_{0}^{K}(\mathbb{F}{})$ as a quotient of $\Mot_{0}%
^{K}(\mathbb{F}{}_{p},n)$ defines an equivalence from $\Mot_{0}^{K}%
(\mathbb{F}{})$ to the tannakian category whose objects are pairs $(V,F)$ as
before together with an action $L\overset{\text{{\tiny def}}}{=}\omega
(X_{mn})$ such that%
\[
F(av)=\sigma a\cdot Fv,\quad a\in L\text{, }v\in V
\]
(cf. \cite{milne2005}, 2.3 and 2.12 et seq.).

\subsection{The category $\Mot^{K}(\mathbb{F}{}_{p},n)$}

Let $F$ be a quadratic extension of $\mathbb{Q}{}$ such that the local degrees
at $p$ and $\infty$ are both $2$. Then $\Mot^{K}(\mathbb{F}{}_{p},n)$ has an
explicit description as an $F$-linear category with a descent datum.

\subsection{The category $\Mot^{K}(\mathbb{F}{})$.}

Again, realize $\Mot^{K}(\mathbb{F}{})$ as a quotient of $\Mot^{K}%
(\mathbb{F}{}_{p},n)$.

\section{Fibre functors on $\Mot_{0}(\mathbb{F}{}_{p})$.}

If each of $\Mot_{0}^{K}(\mathbb{F}{}_{p},n)$ is neutral, does this imply that
$\Mot_{0}(\mathbb{F}{}_{p})=\bigcup_{K,n}\Mot_{0}^{K}(\mathbb{F}{}_{p},n)$ is
neutral? Let $\omega$ be a $\mathbb{Q}{}$-valued fibre functor on
$\Mot_{0}(\mathbb{F}{}_{p})$. Then $\omega$ restricts to a $\mathbb{Q}{}%
$-valued fibre on $\Mot_{0}^{K}(\mathbb{F}{}_{p},n)$ for each $K,n$.

\bquote(Kontsevich email, May 7, 2006). The tower structure means that we have
an epimorphism%
\[
\hat{\mathbb{Z}}^{\times}\rightarrow\Gal(\mathbb{Q}{}^{\mathrm{ab}}%
/\mathbb{Q}{})\rightarrow\hat{\mathbb{Z}}^{\times}\text{.}%
\]

One cannot get $\mathbb{Z}{}_{p}$ factor in the image if one uses only
unramified at $p$ extensions; also if one ignores $\mathbb{Z}{}_{p}$ component
there will be still something wrong: we should get an epimorphism%
\[
\tstyle\prod_{l\neq p}\mathbb{Z}{}_{l}^{\times}/(\mathbb{Z}{}/(l-1)\mathbb{Z}%
{})\twoheadrightarrow\tstyle\prod_{l\neq p}\mathbb{Z}{}_{l}%
\]
which splits the inclusion of the closure of the subgroup generated by the
element $p$. There is a well-known conjecture, 100\% solid by probabilistic
reasons, that for any prime $p$ there are infinitely many primes $l$ such that
$p^{l-1}=1$ mod $l^{2}$, hence $p$ generates a proper closed subgroup in
$\mathbb{Z}{}_{l}^{\times}/(\mathbb{Z}{}/(l-1)\mathbb{Z}{}=\mathbb{Z}{}_{l}$
by the logarithmic map.\equote

We look at this more generally. Let $\mathsf{M}$ be a tannakian category over
$k$ that is a countable union $\mathsf{M}=\bigcup\mathsf{M}_{n}$,
$\mathsf{M}_{n}\subset\mathsf{M}_{n+1}$, of neutral algebraic tannakian subcategories.

Suppose first that $k$ is algebraically closed, and chose a $k$-valued fibre
functor $\omega_{n}$ on each $\mathsf{M}_{n}$. Because $k$ is algebraically
closed, $\omega_{n+1}|\mathsf{M}_{n}\approx\omega_{n}$. In fact, given
$\omega_{n}$, we can modify $\omega_{n+1}$ so that $\omega_{n+1}%
|\mathsf{M}_{n}=\omega_{n}$. Thus, there exists a fibre functor $\omega$ on
$\mathsf{M}$ such that $\omega|\mathsf{M}_{n}=\omega_{n}$.

When we try to do this with $k$ not algebraically closed, then we obtain a
sequence of torsors $\underline{\Hom}^{\otimes}(\omega_{n},\omega
_{n+1}|\mathsf{M}_{n})$. Of course, by making a different choice of fibre
functors, we get a different sequence of torsors, but if, for example, the
fundamental groups $P_{n}$ of the $\mathsf{M}_{n}$ are commutative, then we
get in this way a well-defined element of $\varprojlim\nolimits^{1}%
H^{1}(k,P_{n})$, which is the obstruction to $\mathsf{M}$ being
neutral.\footnote{Recall that for an inverse system $(A_{n},u_{n})$ of abelian
groups indexed by $(\mathbb{N},\leq)$, $\varprojlim A_{n}$ and $\varprojlim
^{1}A_{n}$ are the kernel and cokernel respectively of%
\begin{equation}
(\ldots,a_{n},\ldots)\mapsto(\ldots,a_{n}-u_{n+1}(a_{n+1}),\ldots)\colon
\prod\nolimits_{n}A_{n}\xrightarrow{1-u}\prod\nolimits_{n}A_{n}\text{.}%
\end{equation}
}

\section{A replacement for the Tate conjecture}

Let $\mathbb{A}{}^{p,\infty}$ be the restricted product of the $\mathbb{Q}%
{}_{l}$ for $l\neq p,\infty$, and let $\mathbb{A}{}$ be the product of
$\mathbb{A}{}^{p,\infty}$ with the field of fractions of the ring of Witt
vectors with coefficients in the ground field.

\begin{definition}
Suppose that for each variety $X$ in $\mathcal{S}{}$ and each integer $r$ we
have a $\mathbb{Q}{}$-structure $T^{r}(X)$ on the $\mathbb{A}{}$-module
$\mathcal{T}{}^{r}(X)$ of Tate classes. We call the family $(T^{r}(X))_{X,r}$
a \emph{theory of rational Tate classes on }$\mathcal{S}{}$\emph{ }if

\begin{enumerate}
\item for each variety $X$ in $\mathcal{S}{}$, $T^{\ast}(X)\overset
{\text{{\tiny def}}}{=}\bigoplus_{r}T^{r}(X)$ is a $\mathbb{Q}{}$-subalgebra
of $\mathcal{T}{}^{\ast}(X)$;

\item for every regular map $f\colon X\rightarrow Y$ of abelian varieties,
$f_{\ast}$and $f^{\ast}$ preserve the $\mathbb{Q}{}$-structures;

\item every divisor class on $X$ lies in $T^{1}(X)$.
\end{enumerate}
\end{definition}

The elements of $T^{\ast}(X)$ will then be called the rational Tate classes on
$X$ (for the particular theory).

Now let $\mathcal{S}{}$ be the smallest class satisfying the conditions in the
introduction, and assume there exists a theory of rational Tate classes. Then
we can define categories of motives using the varieties in $\mathcal{S}{}$
with the rational Tate classes as the correspondences, and everything in the
preceding sections holds true. If, moreover, algebraic classes are rational
Tate classes, then there is an exact tensor functor from the category of
motives defined by algebraic classes to the category of motives defined by
rational Tate classes. In particular, a fibre functor on the latter gives rise
to a fibre functor on the former.

\bsmall

{
\bibliographystyle{cbe}
\bibliography{/MData/TeX/bib/refs}

\begin{thebibliography}{}

\bibitem[\protect\astroncite{Albert}{1939}]{albert1939}
{\sc Albert, A.~A.} 1939.
\newblock Structure of {A}lgebras.
\newblock American Mathematical Society Colloquium Publications, vol. 24.
  American Mathematical Society, New York.

\bibitem[\protect\astroncite{Hooley}{1967}]{hooley1967}
{\sc Hooley, C.} 1967.
\newblock On {A}rtin's conjecture.
\newblock {\em J. Reine Angew. Math.} 225:209--220.

\bibitem[\protect\astroncite{Lang}{2002}]{lang2002}
{\sc Lang, S.} 2002.
\newblock Algebra, volume 211 of {\em Graduate Texts in Mathematics}.
\newblock Springer-Verlag, New York.

\bibitem[\protect\astroncite{Langlands and Rapoport}{1987}]{langlandsR1987}
{\sc Langlands, R.~P. and Rapoport, M.} 1987.
\newblock Shimuravariet{\"a}ten und {G}erben.
\newblock {\em J. Reine Angew. Math.} 378:113--220.

\bibitem[\protect\astroncite{Lenstra}{1977}]{lenstra1977}
{\sc Lenstra, Jr., H.~W.} 1977.
\newblock On {A}rtin's conjecture and {E}uclid's algorithm in global fields.
\newblock {\em Invent. Math.} 42:201--224.

\bibitem[\protect\astroncite{Milne}{1986}]{milne1986ar}
{\sc Milne, J.~S.} 1986.
\newblock Arithmetic duality theorems, volume~1 of {\em Perspectives in
  Mathematics}.
\newblock Academic Press Inc., Boston, MA.

\bibitem[\protect\astroncite{Milne}{1994}]{milne1994}
{\sc Milne, J.~S.} 1994.
\newblock Motives over finite fields, pp. 401--459.
\newblock {\em In} Motives (Seattle, WA, 1991), Proc. Sympos. Pure Math. Amer.
  Math. Soc., Providence, RI.

\bibitem[\protect\astroncite{Milne}{1997}]{milneCFT}
{\sc Milne, J.~S.} 1997.
\newblock Class field theory.
\newblock Available at www.jmilne.org.

\bibitem[\protect\astroncite{Milne}{2003}]{milne2003}
{\sc Milne, J.~S.} 2003.
\newblock Gerbes and abelian motives.
\newblock Preprint available at www.jmilne.org/math/ (also
  arXiv:math.AG/0301304).

\bibitem[\protect\astroncite{Milne}{2004}]{milne2004}
{\sc Milne, J.~S.} 2004.
\newblock Periods of abelian varieties.
\newblock {\em Compos. Math.} 140:1149--1175.

\bibitem[\protect\astroncite{Milne}{2005}]{milne2005}
{\sc Milne, J.~S.} 2005.
\newblock Quotients of tannakian categories and rational {T}ate classes.
\newblock Preprint, available at www.jmilne.org/math/; also
  arXiv:math.CT/0508479.

\bibitem[\protect\astroncite{Reiner}{2003}]{reiner2003}
{\sc Reiner, I.} 2003.
\newblock Maximal orders, volume~28 of {\em London Mathematical Society
  Monographs. New Series}.
\newblock The Clarendon Press Oxford University Press, Oxford.

\bibitem[\protect\astroncite{Saavedra~Rivano}{1972}]{saavedra1972}
{\sc Saavedra~Rivano, N.} 1972.
\newblock Cat{\'e}gories {T}annakiennes.
\newblock Springer-Verlag, Berlin.

\bibitem[\protect\astroncite{Wei}{1993}]{wei1993}
{\sc Wei, W.} 1993.
\newblock Weil numbers and generating large field extensions.
\newblock PhD thesis, University of Michigan.

\end{thebibliography}
} \esmall

\end{document}